# Inverse Galois Problem and Significant Methods

Fariba Ranjbar[*], Saeed Ranjbar

[*] School of Mathematics, Statistics and Computer Science, University of Tehran, Tehran, Iran.

*fariba.ranjbar@ut.ac.ir*

**ABSTRACT.** *The inverse problem of Galois Theory was developed in the early 1800's as an approach to understand polynomials and their roots. The inverse Galois problem states whether any finite group can be realized as a Galois group over $\mathbb{Q}$ (field of rational numbers). There has been considerable progress in this as yet unsolved problem. Here, we shall discuss some of the most significant results on this problem. This paper also presents a nice variety of significant methods in connection with the problem such as the Hilbert irreducibility theorem, Noether's problem, and rigidity method and so on.*

## I. *Introduction*

*Galois Theory was developed in the early 1800's as an approach to understand polynomials and their roots. Galois Theory expresses a correspondence between algebraic field extensions and group theory. We are particularly interested in finite algebraic extensions obtained by adding roots of irreducible polynomials to the field of rational numbers. Galois groups give information about such field extensions and thus, information about the roots of the polynomials. One of the most important applications of Galois Theory is solvability of polynomials by radicals.*

*The first important result achieved by Galois was to prove that in general a polynomial of degree 5 or higher is not solvable by radicals. Precisely, **he stated that a polynomial is solvable by radicals if and only if its Galois group is solvable.***

*According to the Fundamental Theorem of Galois Theory, there is a correspondence between a polynomial and its Galois group, but this correspondence is in general very complicated. The inverse Galois problem deals with this complexity. In particular, since it is difficult to consider the case of a general separable polynomial of degree n, for any integer n, the inverse Galois Theory treats the converse question:*

***Is every finite group realizable as the Galois group of a Galois extension of $\mathbb{Q}$?***

*(A) **General existence problem.** Determine whether G occurs as a Galois group over K. In other words, determine whether there exists a Galois extension M/K such that the Galois group Gal (M/K) is isomorphic to G. We call such a Galois extension M a G-extension over K.*

*(B) **Actual construction.** If G is realisable as a Galois group over K, construct explicit polynomials over K having G as a Galois group. More generally, construct a family of polynomials over a K having G as Galois group.*

---

*The classical Inverse Problem of Galois Theory is the existence problem for the field K = Q of rational numbers.*

*It would of course be particularly interesting if the family of polynomials we construct actually gives all G-extensions of K. One obvious way of formulating this is in the form of a parametric or generic polynomial.*

*The next natural question after (B) one may ask is thus:*

*(C) **Construction of generic polynomials.** Given K and G as above, determine whether a generic polynomial exists for G-extensions over K, and if so, find it. This raises a further question:*

*(D) **The Number of Parameters.** What is the smallest possible number of parameters for a generic polynomial for G-extensions over K? (Again, assuming existence.)*

*Remarks. The existence problem (A) has been solved in the affirmative in some cases. On the other hand, for certain fields, not every finite group occurs as a Galois group.*

*(1) If K = ℂ(t), where t is an indeterminate, any finite group G occurs as a Galois group over K. This follows basically from the Riemann Existence Theorem. More generally, the absolute Galois group of the function field K(t) is free pro-finite with infinitely many generators, whenever K is algebraically closed, [5, 15].*

*(2) If K is a 𝔭-adic field, and K(t) a function field over K with indeterminate t, any finite group G occurs as a Galois group over K(t), by the Harbater Existence Theorem [4].*

*The Inverse Galois Problem is particularly significant when K is the field ℚ of rational numbers (or, more generally, an algebraic number field), or a function field in several indeterminates over ℚ (or over an algebraic number field). The big question is then*

***The Regular Inverse Galois Problem.** Is every finite group realisable as the Galois group of a regular extension of ℚ(t)?*

*Whenever we have a Galois extension M/ℚ(t) (regular or not), it is an easy consequence of the Hilbert Irreducibility Theorem (covered in Chapter III below) that there is a 'specialisation' M/ℚ with the same Galois group. Moreover, if M/ℚ(t) is regular, we get such specialised extensions M/K over any Hilbertian field in characteristic 0, in particular over all algebraic number fields. Hence the special interest in the Regular Inverse Galois Problem. Inverse Problem of Galois Theory has been a difficult problem; it is still unsolved.*

**II.** *Milestones in Inverse Galois Theory*

*The Inverse Galois Problem was perhaps known to Galois. In the early nineteenth century, the following result was known as folklore:*

***The Kronecker-Weber Theorem.** Any finite abelian group G occurs as a Galois group over Q: Indeed G is realized as the Galois group of a subfield of the cyclotomic field ℚ(ξ), where ξ is an $n^{th}$ root of unity for some natural number n [9].*

*The first systematic study of the Inverse Galois Problem started with Hilbert in 1892. Hilbert used his Irreducibility Theorem (see Chapter III) to establish the following results [6]:*

**Theorem.** *For any $n \geq 1$, the symmetric group $S_n$ and the alternating group $A_n$ occur as Galois groups over Q.*

*In 1916, E. Noether [13] raised the following question:*

**The Noether Problem.** *Let $M = \mathbb{Q}(t_1, \ldots, t_n)$ be the field of rational functions in n indeterminates. The symmetric group $S_n$ of degree n acts on M by permuting the indeterminates. Let G be a transitive subgroup of $S_n$, and let $K = M^G$ be the subfield of G-invariant rational functions of M. Is K a rational extension of $\mathbb{Q}$? I.e., is K isomorphic to a field of rational functions over $\mathbb{Q}$? (see Chapter III).*

*If the Noether Problem has an affirmative answer, G can be realized as a Galois group over $\mathbb{Q}$, and in fact over any Hilbertian field of characteristic 0, such as an algebraic number field.*

*The next important step was taken in 1937 by A. Scholz and H. Reichardt [18, 16] who proved the following existence result:*

**Theorem 2.** *For an odd prime p, every finite p-group occurs as a Galois group over $\mathbb{Q}$.*

*The final step concerning solvable groups was taken by Shafarevich [21] (with correction appended in 1989; for a full correct proof, the reader is referred to Chapter IX of the book by Neukirch, Schmidt and Wingberg [14]), extending the result of Iwasawa [8] that any solvable group can be realized as a Galois group over the maximal abelian extension $\mathbb{Q}^{ab}$ of $\mathbb{Q}$.*

**Theorem (Shafarevich).** *Every solvable group occurs as a Galois group over $\mathbb{Q}$.*

*Shafarevich's argument, however, is not constructive, and so does not produce a polynomial having a prescribed finite solvable group as a Galois group.*

**Some remarks regarding simple groups.** *Of the finite simple groups, the projective groups PSL(2, p) for some odd primes p were among the first to be realized. The existence was established by Shih in 1974, and later polynomials were constructed over $\mathbb{Q}(t)$ by Malle and Matzat:*

**Theorem [22].** *Let p be an odd prime such that either 2, 3 or 7 is a quadratic non-residue modulo p. Then PSL(2, p) occurs as a Galois group over $\mathbb{Q}$.*

*For the 26 sporadic simple groups, all but possibly one, namely, the Mathieu group $M_{23}$, have been shown to occur as Galois groups over $\mathbb{Q}$. For instance:*

**Theorem (Matzat & al. [12]).** *Four of the Mathieu groups, namely $M_{11}$, $M_{12}$, $M_{22}$ and $M_{24}$, occur as Galois groups over $\mathbb{Q}$.*

*Matzat and his collaborators further constructed families of polynomials over $\mathbb{Q}(t)$ with Mathieu groups as Galois groups.*

*The most spectacular result is, perhaps, the realization of the Monster group, the largest sporadic simple group, as a Galois group over $\mathbb{Q}$ by Thompson [25]. In 1984, Thompson succeeded in proving the following existence theorem:*

***Theorem (Thompson).*** *The monster group occurs as a Galois group over $\mathbb{Q}$.*

*Most of the aforementioned results dealt with the existence question (A) for $K = \mathbb{Q}$.*

*It should be noted that all these realization results of simple groups were achieved via the rigidity method and the Hilbert Irreducibility Theorem (see Chapter III).*

## III. *Significant Methods*

### 1. *The Hilbert irreducibility theorem*

***Definition.*** *Let K be a field, and let $f(t, x)$ be an irreducible polynomial in $K(t)[x] = K(t_1, ..., t_r)[x_1, ..., x_s]$. We then define the Hilbert f-set $H_f/K$ as the set of tuples $a = (a_1, ..., a_r) \in K^r$ such that $f(a, x) \in K[x]$ is well-defined and irreducible. Furthermore, we define a Hilbert set of $K^r$ to be the intersection of finitely many Hilbert f-sets and finitely many subsets of $K^r$ of the form $\{a | g(a) \neq 0\}$ for a nonzero $g(t) \in K[t]$.*

*The field K is called* **Hilbertian**, *if the Hilbert sets of $K^r$ are non-empty for all r. In this case, they must necessarily be infinite.*

*Let K be a field of characteristic 0. Then the following conditions are equivalent:*

*(i) K is Hilbertian.*

*(ii) If $f(t, X) \in K[t, X]$ has no roots in $K(t)$ (as a polynomial in X) there is an $a \in K$ such that f(a,X) has no roots in K.*

***The Hilbert irreducibility Theorem.*** *$\mathbb{Q}$ is Hilbertian.*

***Sketch of proof.*** *Let $f(t, X) \in \mathbb{Z}[t, X]$ have degree n in X, and assume that it has no roots in $\mathbb{Q}(t)$. First thing to do is to translate f(t,X) such that 0 becomes a regular point. Then we have root functions $\theta_1, ..., \theta_n$ defined on a neighborhood of $\infty$, i.e., for all t with |t| greater than some $T \in \mathbb{R}_+$ and*

$$f(t, X) = g(t)\Pi_{i=1}^{n}(X - \theta_i(t)) \in \mathbb{C}((t^{-1}))[X],$$

*where $g(t) \in \mathbb{Q}[t]$ is the coefficient of $X^n$ in f(t,X). Now, if $a \in \mathbb{Q}$ (with a > T ) is such that none of $\theta_1(a), ..., \theta_n(a)$ are rational, therefore obviously f(a,X) has no rational roots. We must prove that there is such an a.* □

***Result.*** *Let K be a Hilbertian field. If a finite group G occurs as a Galois group over K(t), it occurs over K as well.*

***Theorem.*** *Let K be a Hilbertian field, and let $f(t, X) \in K(t)[X]$ be monic, irreducible and separable. Then there is a Hilbert set of $K^r$ on which the specialisations $f(a, X) \in K[X]$ of $f(t, X)$ are well-defined, irreducible and*

$$Gal(f(a, X)/K) \simeq Gal(f(t, X)/K).$$

*By these results, Hilbert proved that for any positive integer n, the symmetric group $S_n$ and the alternating group $A_n$ are Galois groups over $\mathbb{Q}$.*

## 2. Noether's problem

***Definition.*** *If an extension is not algebraic, then we call it, transcendental.*

***Theorem (Emmy Noether).*** *If G is finite and $\mathbb{Q}(X)^G/\mathbb{Q}$ is rational (purely transcendental), then there is a Galois field extension K/$\mathbb{Q}$ with group G.*

***Proof.*** *If $\mathbb{Q}(X)^G/\mathbb{Q}$ is rational, then $\mathbb{Q}(X)^G \cong \mathbb{Q}(w_1, \ldots, w_n)$ where the $w_i$'s are indeterminates. Notice that $\mathbb{Q}(X)/(\mathbb{Q}(X))^G$ is Galois with group G. Thus, there is an irreducible polynomial $f(y) \in \mathbb{Q}(w_1, \ldots, w_n)$ such that $\mathbb{Q}(X)$ is the splitting field for f(y). Using the Primitive Element Theorem, there is an $\alpha \in \mathbb{Q}(X)$ such that $f(\alpha) = 0$, and $\mathbb{Q}(X) = \mathbb{Q}(X)^G(\alpha)$.*

*Consider $t_0$ be a point in $\mathbb{Q}^n$. Let $f_{t_0}(y)$ be the polynomial in $\mathbb{Q}[y]$ obtained by replacing $t_0$ by $(w_1, \ldots, w_n)$ in $f(y) \in \mathbb{Q}(w_1, \ldots, w_n)[y]$. Using results of Hilbert's Irreducibility Theorem, we can have infinitely many $t_0 \in \mathbb{Q}^n$ such that $f_{t_0}(y) \in \mathbb{Q}[y]$ is irreducible (over $\mathbb{Q}$). Now assume that $L \supset \mathbb{Q}$ be the splitting field of $f_{t_0}(y)$ for some $t_0$ where $f_{t_0}(y)$ is irreducible. Therefore L is Galois over $\mathbb{Q}$ with group G.* □

## 3. The embedding problem

*Let $1 \to A \to \tilde{G} \to G \to 1$ be an exact sequence of finite groups. Suppose that G is a Galois group over a field K, $G \cong Gal(L/K)$. The question is whether there is a Galois extension $\tilde{L}$ of K such that $\tilde{G} \cong Gal(\tilde{L}/K)$, $\tilde{L} \supset L \supset K$ and the diagram*

$$\begin{array}{ccccccccc} 1 & \to & A & \to & \tilde{G} & \to & G & \to & 1 \\ & & \downarrow & & \downarrow & & \downarrow & & \\ 1 & \to & Gal(\tilde{L}/L) & \to & Gal(\tilde{L}/K) & \to & Gal(L/K) & \to & 1 \end{array}$$

*is commutative [7].*

*The group A is called the kernel of the embedding problem.*

## 4. Rigidity method (Rationality criteria)

*A great progress has been made in the realization of simple groups as Galois group of regular extensions over $\mathbb{Q}(t)$ and, by Hilbert's irreducibility theorem, over every number field. The rigidity method validates a classical idea: Using the fact that every finite group is a Galois group of a polynomial with coefficient in $\mathbb{C}(t)$ and impose conditions in order to ensure that the polynomial can be defined over $\mathbb{Q}(t)$.*

*First we consider $\chi = \hat{\mathbb{C}}$ or equivalently $\mathbb{P}^1$ and let $S = \{p_1, \ldots, p_s\} \subset \mathbb{P}^1$ be a finite set of removed points. Then the fundamental group $\pi_1(\mathbb{P}^1 \backslash S; p_0)$ with respect to any base point $p_0 \in \chi \setminus S$ has the structure*

$$\pi_1 = \pi_1(\mathbb{P}^1 \backslash S; p_0) = \langle \gamma_1, \ldots, \gamma_s | \gamma_1 \ldots \gamma_s = 1 \rangle.$$

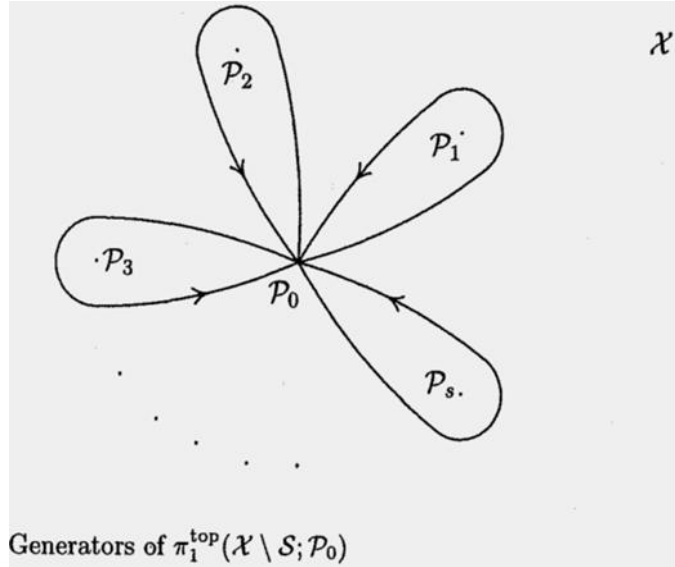

Generators of $\pi_1^{top}(\mathcal{X} \setminus \mathcal{S}; \mathcal{P}_0)$

Then if $\mathbb{C}(t)_S$ be the maximal Galois extension of $\mathbb{C}(t)$ unramified outside $S$, we have:

$$G_S = Gal(\mathbb{C}(t)_S/\mathbb{C}(t)) \cong \hat{\pi}_1 \text{ (\textbf{The existence Riemann theorem})},$$

Where $\hat{\pi}_1$ is the profinite completation of $\pi_1$. Therefore $G_S$ is a profinite group with s generators and the relation $\gamma_1 \ldots \gamma_s = 1$. Let G be a finite group. We can always consider $g_1, \ldots, g_s$ generators with relation $g_1 \ldots g_s = 1$. Then we can have the following epimorphism:

$$\psi: G_S \longrightarrow G, \psi(\gamma_i) = g_i.$$

Now let $N = (\mathbb{C}(t)_S)^{ker(\psi)}$ be fixed field, we have

$$Gal(\mathbb{C}(t)_S/\mathbb{C}(t))/ker(\psi) \cong Gal(N/\mathbb{C}(t)) \cong G.$$

Notice that extension $N/\mathbb{C}(t)$ is unramified outside S because it is obviously $N \subset \mathbb{C}(t)_S$. Thus we have:

**Corollary [10].** *Every finite group is the Galois group of some Galois extension of $\mathbb{C}(t)$ and the Galois group of some Galois extension of $\mathbb{R}(t)$.*

**Theorem [10].** *Every finite group is the Galois group of some extension of $\overline{\mathbb{Q}}(t)$, where $\overline{\mathbb{Q}}$ is the algebraic closure of $\mathbb{Q}$.*

*The difficult problem is to descend to $\mathbb{Q}$ which will be achieved by rigidity method.*

**Definition (Rigid).** *Let G be a finite group. Let $C_1, \ldots, C_r$, $r \geq 3$, be a r-tuple of conjugacy classes of G. Let us denote*

$$\bar{A} = \bar{A}(C_1, \ldots, C_r) = \{(g_1, \ldots, g_r) \in C_1 \times \ldots \times C_r : g_1 \ldots g_r = 1\},$$

$$A = A(C_1, \ldots, C_r) = \{(g_1, \ldots, g_r) \in \bar{A} | \langle g_1, \ldots, g_r \rangle = G\},$$

$$\mathbb{Q}_C = \mathbb{Q}_{C_1, \ldots, C_r} = \mathbb{Q}_{x_1, \ldots, x_r} = \mathbb{Q}(\{\chi(C_i) | \chi \in Irr(G), i = 1, \ldots, r\})$$

*where $x_i \in C_i$, $i = 1, \ldots, r$ and Irr(G) denotes the set of irreducible character of G. Clearly $A \subset \bar{A}$ and G operates by conjugacy on A and on $\bar{A}$.*

The family $(C_1, \ldots, C_r)$ (or $A = A(C_1, \ldots, C_r)$) is called rigid if $A$ is not empty and $G$ acts transitively on $A$.

The family $(C_1, \ldots, C_r)$ (or $A = A(C_1, \ldots, C_r)$) is called strictly rigid if it is rigid and $A = \bar{A}$.

$A$ is rationally rigid if it is rigid and $C_i$ is rational for $i = 1, \ldots, r$ (a conjugacy class $C$ of $G$ is called rational over $\mathbb{Q}$ if any irreducible character of $G$ is rational on $C$).

**Lemma [2].** Suppose that the center of $G$ is $\{1\}$. Then $A = A(C_1, \ldots, C_r)$ is rigid if and only if $|G| = |A|$; and $A = A(C_1, \ldots, C_r)$ is strictly rigid if and only if $|G| = |\bar{A}|$.

The cardinality of $\bar{A} = \bar{A}(C_1, \ldots, C_r)$ can be computed if the character table of $G$ is known, so we have

$$|\bar{A}| = \frac{|C_1| \ldots |C_r|}{|G|} \left( \sum_{\chi \in Irr(G)} \frac{\chi(x_1) \ldots \chi(x_r)}{\chi(1)^{r-2}} \right)$$

where $x_i \in C_i$, $i = 1, \ldots, r$ and $Irr(G)$ denotes the set of irreducible character of $G$ [20]. The rationality of the conjugacy classes can be also checked from the character table of $G$.

Suppose that a finite group $G$ with center $\{1\}$, has a family $(C_1, \ldots, C_r)$ rationally rigid. Let $\psi$ be the mentioned epimorphism. Let $N = (\bar{\mathbb{Q}}(t)_S)^{ker(\psi)}$, $N$ is a Galois extension of $\bar{\mathbb{Q}}(t)$ with $Gal(N/\bar{\mathbb{Q}}(t)) \cong G$. The purpose of rational and rigid conditions is that they imply that $N/\mathbb{Q}(t)$ is normal and the Galois group $\Gamma = Gal(N/\mathbb{Q}(t))$ contains a subgroup $H \subset \Gamma$ such that $\Gamma = HG$. Thus, the fixed field $N_\circ = N^H$ is a Galois extension of $\mathbb{Q}(t)$ with Galois group $G$ such that $N_\circ \bar{\mathbb{Q}} = N$.

**Rigidity Theorem [1, 3, 10, 11, 25].** Let $G$ be a finite group with center $\{1\}$. Let $A = A(C_1, \ldots, C_r)$, $K = \mathbb{Q}_C = \mathbb{Q}_{C_1,\ldots,C_r}$ and let $t$ be an indeterminate over $K$. Assume that $A$ is rigid. Then for any arbitrary chosen set $S = \{p_1, \ldots, p_r\}$ of prime divisors of $\mathfrak{P}_i \in \mathbb{P}(\mathbb{Q}_C(t)/\mathbb{Q}_C)$ of degree one there exists a Galois extension $N$ of $K(t)$ unramified outside $S$ with

$$Gal(N/\mathbb{Q}_C(t)) \cong G$$

Such that the inertia groups over the $\mathfrak{P}_i$ are generated by elements $\sigma_i \in C_i$. If $A$ is rationally rigid, we have $\mathbb{Q}_C = \mathbb{Q}$.

The Hilbert Irreducibility Theorem now implies

**Corollary.** With assumptions and notations as in rigidity theorem, $G$ is a Galois group over $K$.

### 5. Arithmetic – Geometric method (Elliptic curves and the group $GL_2(\mathbb{F}_p)$)

**Definition.** An **elliptic curve** is a smooth, projective algebraic curve of genus one, on which there is a specified point $O$. An elliptic curve is in fact an abelian variety that is, it has a multiplication defined algebraically, with respect to which it is a (necessarily commutative) group and $O$ serves as the identity element.

*Any elliptic curve can be written as a plane algebraic curve defined by an equation of the form:* **$y^2 = x^3 + ax^2 + bx + c$** *which is non-singular; that is, its graph has no cusps or self-intersections. The point O is actually the **"point at infinity"** in the projective plane.*

### *The Mordell group law:*

*If P and Q are two points on the curve (i.e. their coordinates are solutions to the equation of the curve), then we can uniquely describe a third point which lies in the intersection of the curve with the line through P and Q.*

*It is then possible to introduce a group operation, +, on the curve with the following properties: we consider the point at infinity to be 0, the identity of the group; and if a straight line intersects the curve at the points P, Q and R, then we require that P + Q + R = 0 in the group. One can check that this turns the curve into an abelian group, and thus into an abelian variety [23].*

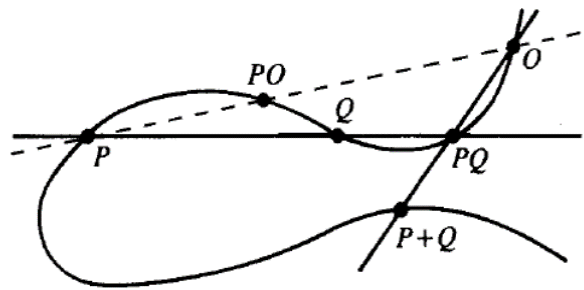

*Mordell group law applied to points of an elliptic curve: this means P+Q is the third intersection point on the line through 0 and PQ.*

**Definition (K-points (Rational points)).** *Let K be a field over which the curve is defined (i.e., the coefficients of the defining equation or equations of the curve are in K) and denote the curve by E. Then the K-points of E are the points on E whose coordinates all lie in K, including the point at infinity. The set of K-points is denoted by E(K). It also forms a group.*

*We fix a positive integer $n \in N$. We want to study the set of points of $E(\mathbb{C})$ that are of order dividing n . Denote by [n] the multiplication-by-n map that sends P to [n]P, where [n]P denotes the sum of n times P. Note that the map [n] is well-defined over $E(\mathbb{C})$ by the Mordell group law.*

**Definition (Points of order dividing n).** *Let $P \in E(\mathbb{C})$. We say that P is a point of order dividing n if $P \in Ker([n])$, i.e. if $[n]P = 0$. We denote E[n] the set of points of order dividing n, i.e. $Ker([n]) = E[n]$.*

**Proposition [24].** *The group E[n] is isomorphic to the direct sum of two cyclic groups of order n:*

$$E[n] \cong \mathbb{Z}/n\mathbb{Z} \oplus \mathbb{Z}/n\mathbb{Z}.$$

*In particular, E[n] is generated by two points $P_1, P_2$ and that for any $P \in E[n]$ there exists $a_1, a_2 \in \mathbb{Z}/n\mathbb{Z}$ such that $P = a_1 P_1 + a_2 P_2$. Precisely, $P_1$ and $P_2$ are the images of some generators of $\mathbb{Z}/n\mathbb{Z} \oplus \mathbb{Z}/n\mathbb{Z}$.*

***Proposition [26].*** *Let E be an elliptic curve defined by the Weierstrass equation $y^2 = x^3 + ax^2 + bx + c$, with $a, b, c \in \mathbb{Q}$. Then, write $E[n] = \{O, (x_1, y_1), \ldots, (x_m, y_m)\}$, with $m = n^2 - 1$. Let*

$$K = \mathbb{Q}(E[n]) = \mathbb{Q}(x_1, y_1, x_2, \ldots, x_m, y_m)$$

*be the extension of $\mathbb{Q}$ generated by the coordinates $x_i$ and $y_i$ with $i \in \{1, \ldots, m\}$. Then $K/\mathbb{Q}$ is a Galois extension of $\mathbb{Q}$.*

***Theorem [23].*** *Fix an integer $n \geq 2$. Let E be an elliptic curve given by the Weierstrass equation $y^2 = x^3 + ax^2 + bx + c$, with $a, b, c \in \mathbb{Q}$. Let $P_1$ and $P_2$ be two generators for $E[n]$. Then there exists a one-to-one group homomorphism*

$$\rho_n: Gal(\mathbb{Q}(E[n])/\mathbb{Q}) \longrightarrow GL_2(\mathbb{Z}/n\mathbb{Z}).$$

***Sketch of proof.*** *For proof, we construct the group homomorphism $\rho_n$. For any $\sigma \in Gal(\mathbb{Q}(E[n])/\mathbb{Q})$ we can define a group homomorphism $\Psi(\sigma)$ on $E[n]$. The map $\Psi: Gal(\mathbb{Q}(E[n])/\mathbb{Q}) \longrightarrow Aut(E[n])$ is a group homomorphism. Furthermore, $\sigma$ can be expressed by a matrix:*

$$\Phi: Aut(E[n]) \longrightarrow GL_2(\mathbb{Z}/n\mathbb{Z})$$
$$\sigma \mapsto \Phi(\sigma) = \begin{pmatrix} \alpha & \beta \\ \gamma & \delta \end{pmatrix},$$

*where $\alpha, \beta, \gamma, \delta \in \mathbb{Z}/n\mathbb{Z}$ are given by the definition of $\sigma(P_1) = \alpha P_1 + \gamma P_2$ and $\sigma(P_2) = \beta P_1 + \delta P_2$ (any homomorphism $\sigma: E[n] \longrightarrow E[n]$ is defined by the images of $P_1$ and $P_2$. In particular we know that $\sigma(P_1), \sigma(P_2)$ are in $E[n]$). We easily can see that $\Phi$ is a group isomorphism. Finally, by composing these two homomorphisms, we obtain the following group homomorphism and all we need now is to show that it is one-to-one:*

$$\rho_n: Gal(\mathbb{Q}(E[n])/\mathbb{Q}) \longrightarrow GL_2(\mathbb{Z}/n\mathbb{Z})$$
$$\sigma \mapsto \Phi \circ \Psi(\sigma). \square$$

***Definition (Galois representation).*** *The map*

$$\rho_n: Gal(\mathbb{Q}(E[n])/\mathbb{Q}) \longrightarrow GL_2(\mathbb{Z}/n\mathbb{Z})$$

*defined as in the previous result is called the Galois representation.*

***Proposition [17].*** *Let E be an elliptic curve defined over $\mathbb{Q}$ by the equation $y^2 = x^3 + ax^2 + b$. Let $p \neq 2$ be a prime number and assume that the representation $\rho_p: Gal(\mathbb{Q}(E[p])/\mathbb{Q}) \longrightarrow GL_2(\mathbb{Z}/p\mathbb{Z})$ is surjective, then we have:*

*Let $P = (x, y) \in E[p] - \{0\}$. The characteristic polynomial of the multiplication by $x + y$ in $\mathbb{Q}(x, y)$ is irreducible and its Galois group over $\mathbb{Q}$ is $GL_2(\mathbb{Z}/p\mathbb{Z})$.*

*Serre [19] proved the following theorem which says that under which condition the map $\rho_n$ is onto for almost every $n \in \mathbb{N}$. The proof was presented in 1972 in and need more advanced tools of elliptic curves theory.*

***Theorem.*** *Let E be an elliptic curve given by a Weierstrass equation with rational coefficients. Assume that E does not have complex multiplication. There is an integer*

$N_E \geq 1$ *depending on the elliptic curve E, such that if n is an integer relatively prime to $N_E$, then the Galois representation*

$$\rho_n: Gal(\mathbb{Q}(E[n])/\mathbb{Q}) \longrightarrow GL_2(\mathbb{Z}/n\mathbb{Z})$$

*is an isomorphism.*

## *References*

[1] G.V. Belyi, On Galois extensions of a maximal cyclotomic field, Math. USSR Izvestia A.M.S. Translation 14 (1979) 247-256.

[2] W. Feit, Rigidity and Galois groups, Proc. of the Rutgers group theory year 1983-1984, Cambridge University Press (1984), 283-287.

[3] M. D. Fried, Fields of Definition of Function Fields and Hurwitz Families–Groups as Galois groups, Comm. Algebra **5** (1977), 17-82.

[4] D. Harbater, Galois coverings of the arithmetic line, Number Theory: New York, 1984–85, Lecture Notes in Mathematics. **1240**, Springer-Verlag (1987), 165–195.

[5] D. Harbater, Fundamental groups and embedding problems in characteristic p, Recent Developments in the Inverse Galois Problem (Seattle, WA, 1993), Cotemp. Math. **186** (1995), 353–369.

[6] D. Hilbert, Ueber die Irreduzibilität ganzer rationaler Funktionen mit ganzzahligen Koeffizienten, J. reine angew. Math. **110** (1892), 104–129.

[7] M. Ikeda, Zur Existenz eigentlicher galoisscher Körper Beim Einbettungsproblem für galoissche Algebren, Abh. Math. Sern. Univ. Hamburg **24** (1960), 126-131.

[8] K. Iwasawa, On solvable extensions of algebraic number fields, Ann. Math. 58 (1953), 126–131.

[9] F. Lorenz, Algebraische Zahlentheorie, B. I. Wissenschaftsverlag, Mannheim (1993).

[10] G. Malle, and B. H. Matzat, Inverse Galois Theory, Berlin, Heidelberg, New York, Springer-Verlag (1999).

[11] B. H. Matzat, Konstruktion von Zahl-und Funktionenkörpern mit vorgegebener Galoisgruppe, J. reine angew. Math. **349** (1984), 179-220.

[12] B. H. Matzat, and A. Zeh-Marschke, Realisierung der Mathieugruppen $M_{11}$ und $M_{12}$ als Galoisgruppen über $\mathbb{Q}$, J. Number Theory **23** (1986), 195-202.

[13] E. Noether, Gleichungen mit Vorgeschriebener Gruppe, math. Ann. **78** (1916), 221-229.

[14] J. Neukirch, A. Schmidt, and K. Wingberg, Cohomology of number fields, Grundlehren der Mathematischen Wissenschaften 323, Springer-Verlag (2000).

[15] F. Pop, Étale coverings of affine smooth curves. The geometric case of a conjecture of Shafarevich. On Abhyankar's conjecture, Invent. Math. **120** (1995), 555–578.

[16] H. Reichardt, Konstruktion von Zahlkörpern mit gegebener Galoisgruppe von Primzahlpotenzordnung, J. reine angew. Math. 177 (1937), 1-5.

[17] A. Reverter, and N. Vila, Polynomials of Galois representations attached to elliptic curves, Rev. R. Acad. Cienc. Exact. Fis. Nat. (Esp), Vol. 94. 3º (2000), 417-421.

[18] A. Scholz, Konstruktion algebraischer Zahlkörper mit beliebiger Gruppe von Primzahlpotenzordnung I, Math. Z. **42** (1937), 161-188.

[19] J. P. Serre, Propriétés galoisiennes des points d'ordre fini sur des courbes elliptiques, Berlin, Springer-Verlag, Inventiones math. **15** (1972), 259-331.

[20] J. P. Serre, Topics in Galois Theory, USA, AK Peters Ltd, (1992).

[21] I. R. Shafarevich, Construction of fields of algebraic numbers with given solvable Galois group (in Russian), Izv. Akad. Nauk SSSR, Ser. Mat. **18** (1954), 525-578.

[22] K. _Y. Shih, On the construction of Galois extensions of function fields and number fields, Math. Ann. **207** (1974), 99-120.

[23] J. Silverman, The Arithmetic of Elliptic Curves , New York, Springer, (1986).

[24] J. Silverman and J. Tate, Rational points on elliptic curves , New York, Springer, (1992).

[25] J. G. Thompson, Some finite groups which appear as Gal(L/K), where K ⊆ $\mathbb{Q}(\mu_n)$, J. Algebra **89** (1984), 437-499.

[26] L. Washington, Elliptic: number theory and cryptography, Boca Raton FL, Chapman &Hall (2008).